\documentclass[12pt]{article}
\usepackage{epsfig,amsmath,latexsym}
\usepackage{amsfonts,amssymb,graphicx}
\usepackage{amscd, amssymb, amsthm, amsmath,eucal, verbatim}
\usepackage{epstopdf}
\newtheorem{theorem}{Theorem}[section]

\newtheorem{lemma}[theorem]{Lemma}
\newtheorem{proposition}[theorem]{Proposition}
\newtheorem{corollary}[theorem]{Corollary}

\newcommand{\bproof}{\noindent {\bf Proof. }}

\newcommand{\RR}{\mathbb{R}}

\newcommand{\RP}{\mathbb{RP}}

\begin{document}
\title{Knots \& $k$-width}
\author{Joel Hass \footnote {Supported in part by NSF grants.},
J. Hyam Rubinstein \footnote{Supported in part by the Australian Research Oouncil} and Abigail Thompson \footnote{Supported in part by NSF grants.}.  }
\maketitle

\begin{abstract}
We investigate several natural integer invariants of curves  in $\RR^3$ and explore their values on isotopy classes of curves.
\end{abstract}

\section{Introduction} 
In this paper we introduce the notion of $k$-width for a knot or link in $\RR^3$, where $k$ is an integer between 1 and 4.  These widths provide increasingly detailed information, as $k$ increases, on the
intersections of a curve with flat planes and round spheres in $\RR^3$.
We examine properties of curves that minimize $k$-width within their isotopy class.

The notion of $k$-width is motivated by considerations in the theory of
index-$k$ minimal submanifolds. Our notion of $k$-width  does not appear to have been studied before for $k>1$. 
After introducing $k$-width, we explore properties of $2-$width in some detail. 
We show that only finitely many knots have $2-$width
less than $c$, for any positive constant $c$.  Thus, like crossing number, $2-$width gives a way to order all knots in terms of increasing complexity. We also show that only the unknot and the trefoil have $2-$width less than $10$. 

We also relate the 2-width of a curve to its curvature.  
Milnor and Fary showed that if a smooth curve $\gamma$ has total curvature less than $4\pi$ then $\gamma$ is unknotted  \cite{Milnor} \cite{Fary:49}.  The same argument shows that if the total curvature of $\gamma$ is less than $2n \pi$ then there is a direction relative to which $\gamma$ has at most $n$ maxima.  As a consequence, its bridge number is at most $n$.  The converse is false.  A curve can be deformed to have arbitrarily large curvature without changing the bridge number.  We show here that if an immersed plane curve $\gamma$ has $2-$width $n$ then it has  total curvature at most $2\pi n^{3/2} $.
We prove that if  a curve $\gamma$ in $\RR^3$
 has $2-$width $n$, then some planar projection of $\gamma$ 
 has  total curvature at most $2\pi n^{3/2} $.
This can be viewed in contrast to the Fary-Milnor Theorem.  While small bridge number does not imply that some projection has small total curvature, small 2-thickness does imply this.

In Section~\ref{definitions} we introduce $k$-width for curves in $\RR^3$. In Section~\ref{bounded} we show that there are only a finite number of knot types with 2-width bounded by a given constant and in Section~\ref{low}  we look at the knots with the lowest 2-width.  In Section~\ref{curvature}  we establish a connection between 2-width and the total curvature of a curve. Small 2-width implies that the total curvature of some projection of a curve is small. Finally in Section~\ref{general} we discuss a far more general extension of the notion of width, involving the intersection between a submanifold of a manifold and a group of submanifolds.  

\section{Width for curves in $\RR^3$}  \label{definitions}
  
Combinatorial knot  invariants computed by counting intersections of curves with planes originate with Schubert's bridge number \cite{Schubert}, and also include
Kuiper's superbridge number \cite{Kuiper} and Gabai's thin position \cite{Gabai}.
Let $\gamma$ be a smooth curve in $\RR^3$.
The 1-width of $\gamma$ is found in \cite{Gabai}. It is computed by counting the sum of a curve's intersections with a finite collection of horizontal planes in $\RR^3$, one plane being chosen between each pair of horizontal tangents of  $\gamma $.
We present natural definitions of 
 2, 3 and 4-parameter families of planes and spheres that give interesting width invariants. We will describe these in $\RR^3$, but they are equally natural in ${\mathbb S}^3$.

\noindent
{\bf Definition.}
Consider the set $X_2$ of planes in $\RR^3$ perpendicular to the $xy$-plane.
This set is in 1-1 correspondence with the set of straight lines in the plane, which 
is diffeomorphic to  a punctured $\RR P^2$, or equivalently to an open Mobius band.
Assume that a smooth curve $\gamma$ has been perturbed slightly, so that it has no tangent lines parallel to the $z$-axis and its projection onto the $xy$-plane is in general position, with finitely many transverse double points. Moreover we can assume that no plane in $X_2$ is tangent to $\gamma$ at more than two points. 
Set two planes  in $X_2$ equivalent if they
are connected by a path in $X_2$ consisting of planes transverse to $\gamma$.
The {\em width} of a plane transverse to $\gamma$ is the number 
 of intersection points of the plane with $\gamma$. 
Equivalent planes intersect $\gamma$ in the same number of points.
Summing over representatives of each equivalence class of planes transverse to $\gamma$ 
gives an integer $w_2(\gamma)$ called the {\em 2-width} of $\gamma$.  
A curve $\gamma$ is {\em 2-width minimizing} if it minimizes 2-width in its isotopy class. In the next sections we will explore some properties of 2-width.  The 2-width $w_2(K)$ of a knot $K$ is the 
2-width of a 2-width minimizing representative of the knot.

There is a unique plane in $X_2$ which is tangent at every point on the curve $\gamma$, by our assumption that no tangent line is parallel to the $z$-axis. It follows that the set of non-transverse planes forms a singular curve $\lambda$ in the open Mobius band $X_2$. Double tangencies of $\gamma$ with planes  in $X_2$  give double points of $\lambda$ and inflection points of the planar projection of $\gamma$ give rise to cusps on $\lambda$.  Hence the 2-width is finite, since the number of regions in $X_2$ of planes transverse to $\gamma$ is the same as the finite number of complementary regions of $\lambda$ in the open Mobius band $X_2$.
 
It is possible to define 2-width working entirely with the  projection $\pi(\gamma)$ of $\gamma$ onto the $xy$-plane, since changes in its $z$-coordinate do not affect the intersections  of $\gamma$ with a plane in $X_2$.
To define 3-width and 4-width however, we are obliged to work with 3-dimensional representations of a knot.   

\noindent
{\bf Definition.}
Let $X_3$ be the set of all planes in $\RR^3$. $X_3$ is diffeomorphic to a once punctured $\RP^3$. Two planes in $X_3$ are said to be equivalent if they are connected through a path of planes transverse to $\gamma $.  
The {\em width} of a plane transverse to $\gamma$ is the number 
of intersection points of the plane with $\gamma$. The {\em 3-width} $w_3(\gamma)$ of $\gamma $ 
is the sum  of the widths over all equivalence classes of planes transverse to $\gamma$. 
A curve is {\em 3-width minimizing} if it minimizes 3-width in its isotopy class.

At any point $x$ of $\gamma$, there is a circle of planes touching $\gamma$ at $x$. 
So the set of planes that are not transverse to $\gamma$ forms a singular torus in $X_3$. 
Double curves of this torus are double tangencies of $\gamma$ and triple points are triple tangencies. It is then clear that 3-width is finite, as it is a sum of intersection numbers of planes, one for each of the complementary regions of this singular torus.
 
A related but distinct invariant is the superbridge index introduced by N. Kuiper \cite{Kuiper}.  In our point of view,  bridge number is a 1-parameter and superbridge index is a 3-parameter invariant.  It is possible to define a 2-parameter and 4-parameter bridge index as well, though these do not appear to have been investigated. In \cite{Kuiper} it is shown that there are infinitely many knots with superbridge index bounded by 4. In fact, Kuiper showed that  for all odd $q$, an appropriate embedding of the $(2,q)$-torus knot has superbridge index 4.  
These superbridge indices can be realized by positioning the  $(2,q)$ torus knot as satellites of a ``baseball'' curve on the 2-sphere, a curve  that intersects any plane in at most four points. When positioned on a small radius torus around this baseball curve, the $(2,q)$ torus knot intersects any flat plane in at most 8 points.  Nevertheless, the 3-width of these torus knots does not appear to be bounded.  
   
A 3-width minimizing curve gives one approach to finding an ``optimal" imbedding of a knot.  Note however that the 3-width of a curve is preserved by affine maps of $\RR^3$, and these can radically change the shape of a curve.  Similarly the 2-width of a curve is preserved by  affine maps of $\RR^2$
that preserve the set of planes $X_2$.
   
Finally we consider a 4-parameter width.

\noindent
{\bf Definition.}
Let $X_4$ be the set of all flat planes and round 2-spheres in $\RR^3$. The topology on $X_4$ is obtained by considering the space of round spheres on the 3-sphere, identified to $X_4$ via stereographic projection.  $X_4$ is diffeomorphic to a once punctured $\RP^4$. Define the width of a sphere or plane transverse to $\gamma$ to be the number of its intersection points with $\gamma $. Two spheres or planes are equivalent if they can be connected by a path in $X_4$ consisting of planes and spheres with no tangencies to $\gamma$.  The {\em 4-width} of $\gamma$ is the sum 
 of the widths over all equivalence classes of planes and spheres transverse to 
 $\gamma$. 
A curve is {\em 4-width minimizing} if it minimizes 4-width over all smooth curves in its isotopy class.

As before, one can show that the set of non-transverse spheres and planes is a singular 3-manifold, the image of a product of a Mobius band and a circle. Addition of widths for the finitely many complementary regions gives the 4-width. 

Like the knot energy studied in \cite{FreedmanHeWang:94}, $w_4$ is a conformal invariant of $\gamma$. Unlike knot energies however, $w_4$ is integer valued.

\begin{lemma}
If  $\gamma'$ is the image of  $\gamma$ by a conformal map of $\RR^3$, then $w_4(\gamma') = w_4(\gamma) $.
\end{lemma}
\bproof
There is a width preserving 1-1 correspondence between the equivalence classes for $\gamma$ and those for $\gamma'$.
\qed

We can define higher order widths by considering intersections with quartics and other families of surfaces. See Section~\ref{general} for some very general extensions. The main results of this paper concern 2-width. 
 
\section{The number of knots with bounded 2-width} \label{bounded}
 
Knot tables are usually arranged to reflect increasing crossing number.  This is feasible because for any constant $n$, only finitely many knots have crossing number less than $n$.  This property does not hold for other common invariants,  such as unknotting number, tunnel number and bridge number. While 1-width shares this shortcoming, we prove below a finiteness result for the 2-width invariant.
 
\begin{theorem}
For any constant $n$, only finitely many knots in $\RR^3$ have $2-$width less than or equal to $n$.
\end{theorem}

\bproof
Let $\gamma \subset \RR^3$ be a generic embedded curve. We require additional properties to those in the definition of 2-width, where it was assumed that $\gamma$ has no vertical tangencies and the projection to the horizontal $xy$-plane is in general position. For convenience, we will specify later what we mean by generic.

Recall that the set $X_2$ of planes perpendicular to the $xy$-plane in $\RR^3$ is diffeomorphic to a punctured projective plane. We will not distinguish in what follows
between a point of $X_2$ and the corresponding plane in $\RR^3$.
We associate to $\gamma$ a certain graph $\alpha$ in
$X_2$. Using the projection map $\pi : \RR^3 \to \RR^2$, we project $\gamma$ to a planar curve  $\pi(\gamma)$ in the $xy$-plane. To each point in $\pi(\gamma)$ we associate the tangent line at that point. The set of tangent lines gives a graph $\alpha$ in $X_2$, which we call the  {\em graphic} of $\gamma$ in $X_2$.  (Each plane in $X_2$ corresponds to such a line in the $xy$-plane). The graph $\alpha$ is the image of a curve that is immersed in $X_2$ except at a finite number of cusps. The cusps are the images of inflection points of $\pi(\gamma)$. A finite number of planes in $X_2$ are tangent to $\gamma$ at two points. Such a {\em double tangent plane} corresponds to an order four vertex in $\alpha$. Let $v$ be the number of these vertices.   

Let $r$ count the number of regions of the complement of $\alpha$ in $X_2$. These regions consist of $f$ simply-connected faces,  one non-compact annulus and possibly one Mobius band.   Two points in a given region $R_i$ correspond to planes in $X_2$ in the same equivalence class, so we can assign an integer $w_2(R_i)$ to that region, equal to the number of times a representative plane intersects $\gamma$.  
The 2-width of $\gamma$ is the  sum of the integers  $w_2(R_i)$ over all the regions.

The double tangent planes of $\gamma$ are divided into two types. A double tangent plane is said to be {\em exterior} if it does not separate small neighborhoods of the two tangent points and {\em interior} if it does separate them.  Let $t$ be the number of
exterior double tangent planes of $\gamma$ and $s$ the number of interior
double tangent planes.  Let $i$ be the number of inflection points and let $c$ be the number of crossings of $\pi(\gamma)$.
An inflection point corresponds to a point where $\gamma$ is tangent to a vertical 
plane at a point and the plane separates any neighborhood of the point of
tangency. A crossing point of 
 $\pi(\gamma)$ is a point where a vertical line meets  $\gamma $ twice.
 
A theorem of Fabricius-Bjerre \cite{Fabricius-Bjerre} provides the tool needed to link
properties of $\gamma$ to the graph  $\alpha$. 
Fabricius-Bjerre considered the relation between the numbers of double tangent lines,
crossing points and inflection points of a plane curve $\delta$.  He established the following result \cite{Fabricius-Bjerre}.

\begin{theorem}[Fabricius-Bjerre] \label{Fabricius-Bjerre}
Let $\delta$ be a generic plane curve with $t$ exterior double tangent lines,  $s$ interior
double tangent lines,  $i$  inflection points and  $c$ crossings. Then $c+\frac{1}{2}i=t-s$.
\end{theorem}

In this context, {\em generic} curves are smooth curves in general position, with finite values for each of the quantities in Theorem~\ref{Fabricius-Bjerre}.  For present purposes we can take a generic curve in $\RR^3$ to be one whose projection is generic in this sense.

\begin{lemma}  \label{c<v}
$c \le v$.
\end{lemma}
\bproof
Double tangent planes in $X_2$ correspond to valence-4 vertices of the graph  $\alpha$, so the number of such vertices $v$ of $\alpha$ equals $s+t$. Applying Fabricius-Bjerre's Theorem, we obtain
$$
c \leq c+\frac{1}{2}i=t-s \leq t+s=v .
$$ 
\qed

The next lemma relates $v$ to the number of simply connected faces $f$
in  the complement of $\alpha$. While $\alpha$ may have cusps, they do
not play a role in this computation, so we ignore them.  

\begin{lemma}  \label{f<r<f+3}
$ v = f  $ and either $r=f+1$ or $r = f+2$. The latter occurs when there is a
region in the graphic homeomorphic to a Mobius band.
\end{lemma}
\bproof
There is one region in the complement of $\alpha$ homeomorphic to a half-open annulus. 
There may be one region homeomorphic to a Mobius band.  Annular 
and Mobius band regions do not contribute
to the Euler characteristic computation. Let $e$ count the
number of edges connecting valence four vertices. The Euler characteristic of $X_2$ 
is zero, obtained by taking the sum, giving $ v-e+f = 0$.  Since vertices are of valence four, we have
$e=2v$, and so $f=v$ as claimed.  Depending on whether a Mobius band component exists in the complement of $\alpha$, we have $r=f+1$ or $r = f+2$.
\qed
 
Now let $R_1,R_2,.....R_r$ denote the regions of the graph $\alpha$, with associated widths
$w_1,w_2,...,w_r$ respectively.
Each $w_i$ is a non-negative even integer, and their sum is $n$. At most one of the $w_i$'s equals zero, since there is a single
equivalence class of planes in $\RR^3$ disjoint from $\gamma$.
Thus each contributes at least 2 to the total width of $\gamma$ and the total width $w_2$ satisfies
$w_2(\gamma)  \ge 2r > 2f = 2v \ge 2c $.  

So the crossing number of $\gamma$ is bounded above by half of its $2-$width.  Since only finitely many knots have crossing number less than a given integer, the theorem follows.
\qed

\section{Bounds on 2-width} \label{low}

{\bf Definition.}  A knot projection is {\em positively curved} if it contains
no inflection points.

\noindent
\begin{lemma} \label{positive}
 Every curve in $\RR^3$ can be isotoped to have a positively curved projection.
 \end{lemma}  
\bproof
Every curve $\gamma$ is isotopic to a curve with a braid presentation \cite{Artin}.
Such a presentation is parametrized in cylindrical coordinates by a curve 
$\gamma(s) = (r(s), \theta(s), z(s))$ with $\theta'(s) >0$ and $r(s) >0$.
By scaling we can assume $0 <r(s) < 1$.
An isotopy $\gamma_t(s) = (t + (1-t)r(s), \theta(s), z(s)  ), 0 \le t \le 1- \epsilon$ 
takes $\gamma$ to a curve in a  $2 \epsilon$-neighborhood of the unit radius cylinder $\{r=1\}$.
As $\epsilon \to 0$, $\pi(\gamma)$ smoothly converges to a cover of the unit circle. 
Thus for $\epsilon$ sufficiently small, $\pi(\gamma_{1-\epsilon})$ has all curvatures positive.
\qed

We can apply this to get an upper bound on the 2-width of a knot.

\begin{proposition} 
Let $c$ be the
minimal crossing number of a braid projection of a knot $K$.
Then  $w_2(K) \le  (c+1)(c+2)$, 
\end{proposition}

\noindent
\bproof
Choose $\gamma$ representing $K$ so that $\pi(\gamma$has minimal crossing number over all braid representations of  $K$.
Further isotop  $\gamma$ so that its projection $\pi(\gamma)$ is positively curved as
 in Lemma~\ref{positive}. Then  $\pi(\gamma)$ has no inflection points and no
internal double tangencies.
Therefore its crossing number $c$ equals the number of external double tangencies $t$, which in turn equals the number of vertices $v$ in the graphic of  $\pi(\gamma)$. From Lemma~\ref{f<r<f+3} we know that the number of regions in the graphic $r$ satisfies $ c+1 \le r \le  c+2. $
Each region has width a non-negative even integer, and one
of the regions has width zero.  Regions in the graphic with a common edge have width
that differs by $2$.  So the sum of the widths is at most
$ 2+4+6+...+2(c+1)=(c+2)(c+1). $
\qed

We can obtain a lower bound in terms of the number of intersections between the curve and a 
vertical plane.

\begin{proposition} \label{lb}
If a plane in $X_2$ intersects $\gamma$ in $2n$ 
points then $w_2(\gamma) \ge n(n+1) $
\end{proposition} 
\bproof
Some plane in $X_2$ intersects $\gamma$  in zero points.
Each region in the graphic has width a positive even integer, and
at least one region has width $2n$.  Regions in the graphic
with a common edge have width that differs by $2$, so there are
regions in the graphic of widths $ 2, 4, \dots 2n$ and the total width
is  at least $  2+4+6+...+2n= n(n+1)$.
\qed
 
\begin{theorem} \label{small.width.knots}
The only knots with 2-width less than or equal to 10 are the trefoil
and the unknot.
\end{theorem} 

\bproof
Suppose $w_2(\gamma) \le 10$.
If  $\gamma$ meets a plane in $X_2$ in at least 6 points then
Lemma~\ref{lb} implies that its width is at least 12. So we can
assume that $\gamma$ meets any plane in $X_2$ in at most 4 points.

Call a crossing of  $\pi(\gamma)$ {\em exterior} if it meets the unbounded
region of the complement of $\pi(\gamma)$ in the $xy$-plane. 
Otherwise call it {\em interior}. A ray in the $xy$-plane from a point
near an interior crossing to infinity must cross $\pi(\gamma)$ at least once.  
If $\gamma$ is knotted then $\pi(\gamma)$ has an exterior crossing. Suppose
$\pi(\gamma)$ also has an interior crossing. Then a line in the
 $xy$-plane passing very close to both an interior and an exterior
 crossing meets $\pi(\gamma)$ twice near each crossing and at least
 once more in passing from a point near the interior crossing to infinity.
 Since the intersection number with $\pi(\gamma)$ is even, the intersection
 consists of at least six points. But then $ \gamma$ meets a plane in $X_2$
 in six points and $w_2(\gamma) \ge 12$. So we can assume that all crossings
of  $\pi(\gamma)$ are exterior. 

Color the complementary regions to $\pi(\gamma)$ black and white, 
with the unbounded region white and regions with a common edge
having different colors. Rays based at points in a white region intersect
 $\pi(\gamma)$  in an even number of points, while rays based in
 a black region intersect  $\pi(\gamma)$ an odd number of times.
If there are no bounded white regions then
a small neighborhood of any crossing meets the unbounded region in two components
and a simple induction argument shows that $\gamma$ is unknotted, satisfying the
conclusion of the Theorem. Otherwise there is
at least one bounded white region.
A ray from a bounded white region to infinity meets
$\pi(\gamma)$ in at least two points. If
there are two distinct bounded white regions,  then a line segment connecting a point in each
intersects $\pi(\gamma)$ in at least 2 points.  Extending this
segment to a line by adding two rays gives a line intersecting $\pi(\gamma)$ in at
least six points, contradicting our width assumption. So we can
assume that $\pi(\gamma)$ has exactly one bounded white region, $W$.  

Let $w$ be an interior point of $W$.  Every line through $w$
intersects $\pi(\gamma)$ in exactly 4 points. So $\pi(\gamma)$
gives a $2$-braid projection of $\gamma$, with axis the vertical line
over $w$. The number of crossings of a 2-braid knot is odd,
and since we assume that $\gamma$ is not a trefoil or unknot we
have that the number of crossings is at least 5. 
 Applying Lemma~\ref{c<v} shows that the number
of regions in the graphic is 7, with one annular and one Mobius band region.
One region has width 4 and five other regions have width at least two, giving
that the  2-width of $\gamma$ is at least 14, contradicting our assumption.

We conclude that $\gamma$ is a trefoil or an unknot.
\qed

\begin{theorem}
The $(2,n)$-torus knot $K$ has 2-width $2n+4$. 
\end{theorem} 
\bproof
A positively curved 2-braid projection of $K$ has a graphic with $r = n+2$. 
One region has width 0, one has width 4 and $n$ have width 2, giving a total width of $2n+4$. 
It remains to show no representative of $K$ has smaller 2-width.

Let  $\gamma$ be any representative of $K$. Since $K$ has $n$ crossings in an alternating projection, the number of crossings of $\pi(\gamma)$ is at least $n$.  Applying Lemmas~\ref{c<v} and \ref{f<r<f+3}
shows that the number of regions in the graphic of  $\pi(\gamma)$ is $n+1$ or $n+2$.  As in the proof of
Theorem~\ref{small.width.knots} there is at least one interior white region for $\pi(\gamma)$.  An interior white region has width at least 4.  If there
is exactly one then the graphic has a Mobius band region of width 4 and $r=n+2$. In that case the width is
at least $4 + 2n$ as claimed.  If there is more than one interior white region, then the total width is at least
$4+4+ 2(n-2) = 2n+4$.  The result follows.
\qed

\noindent
{\bf Remark:} Width can be defined and computed for links in the same way. The number of zero-width regions in the graphic associated to a link can be greater than one if the link is split (separated by a plane.) The Hopf link and the 2-component unlink each have $w_2(L) = 8$ and the $(2,4)$-torus link has $w_2(L) = 12$. 

\section{Curvature and 2-width} \label{curvature}

The { \em total curvature} of a curve in $\RR^3$ is obtained by integrating the absolute value of the curvature function. Milnor and Fary showed that if a smooth curve $\gamma$ has total curvature less than $4\pi$ then $\gamma$ is unknotted \cite{Milnor},  \cite{Fary:49}. Milnor's proof proceeds by finding a direction in $\RR^3$ relative to which the curve has only one maximum and one minimum. The same argument shows that if the total curvature of $\gamma$ is less than $2n \pi$ then there is a direction relative to which $\gamma$ has fewer than $n$ maxima.  As a consequence, its bridge number is less than $n$.  The converse is false.  Curves of bridge number two can have unbounded curvature.  In fact any curve can be deformed to have arbitrarily large total curvature without changing its bridge number. So small bridge number does not imply small total curvature. In contrast, we show here that if a curve has small 2-width then its projection onto the associated 2-plane has small total curvature.  

We first obtain a result about immersed plane curves. When a curve lies in the $xy$-plane the notion of 2-width can be restated in terms of its intersections with lines in the plane rather than planes in  $\RR^3$. Lemma~\ref{pos.arc} considers curves immersed in the plane.

\begin{lemma} \label{pos.arc}
Let $\gamma$ be a $C^2$ curve immersed in the plane containing an embedded positively curved arc with total curvature $x$.  Then $w_2(\gamma) >  {x^2 } / ( 2 \pi )^2  .$
\end{lemma} 
\bproof
An embedded positively curved arc with total curvature $x$ is a spiral that winds monotonically around one of its endpoints, with total angle $x$.  For if the arc did not wind monotonically, it would have to touch a ray starting at the given endpoint on one side. But then it is easy to see that the curvature of the arc must have changed sign and this is a contradiction. 

A line through the two endpoints of this spiral intersects  $\gamma$ in at least $[ x / \pi ]+1$ points. Assume first that $[ x / \pi ]+1 $ is even. Then the 2-width of $\gamma $ is at least 
$$
([ x / \pi ]+1 ) +([ x / \pi ]-1 )   + \dots 2  = ([ x / \pi ]+1 )  ([ x / \pi ]+2 ) /4 >  {x^2 } / ( 2 \pi )^2  .
$$

On the other hand, if $[ x / \pi ]+1 $ is odd, then since any line meets $\gamma$ in an even number of points, the line through the endpoints of the spiral intersects $\gamma$ in at least $[ x / \pi ]+2$ points
and the same argument clearly applies. 
\qed

\begin{theorem} \label{planar.curv}
If $\gamma $ is a $C^2$ curve immersed in the plane with total curvature $x$ then $w_2(\gamma) >   Cx^{2/3} $, where $C=1/(2\pi) ^{2/3}$.
\end{theorem} 
\bproof
Let $y = 2c+i$ where $c$ is the number of crossings of  $\gamma$ and $i$ is the number  of its inflection points.  Now by Theorem~\ref{Fabricius-Bjerre} we have that
$2c +i = 2t-2s \le 2t+2s =2v$ where $v$ counts the vertices in the graphic
of $\gamma$. Lemma~\ref{f<r<f+3} implies that 
$ v = f  $ and either $r=f+1$ or $r = f+2$.  So
 $w_2(\gamma) \ge 2f \ge y$. If $y \ge x^{2/3}/(2\pi)^{2/3}$ and the Theorem follows. So
 assume $y < x^{2/3}/  (2\pi)^{2/3}$.
 
The crossing and inflection points of $\gamma $ divide $\gamma $ into $y=2c+i$ subarcs, each with no crossing or inflection points.  Since $\gamma $ has total curvature $x$, one of the subarcs of  
$\gamma$ disjoint from its inflection and crossing points is positively curved with
total curvature greater than $x/y$. By Lemma~\ref{pos.arc} we have that 
$$ 
w_2(\gamma) >   (x/y)^2 /( 2 \pi )^2 =   x^2   / ( 2y \pi )^2  >  
(2\pi)^{4/3} x^2 / ( 2  \pi x^{2/3})^2 =  x^{2/3}/(2\pi)^{2/3} .
$$
\qed

Theorem~\ref{planar.curv} gives the following corollary for curves in $\RR^3$.

\begin{corollary} \label{3d.curv}
If  $\gamma $  is a curve  in $\RR^3$  with $w_2(K) \le  n$ then some planar projection of $ K$ has  total curvature at most $2\pi n^{3/2} $.
\end{corollary} 
\bproof
If a planar projection of  $\gamma $ has curvature greater than $2\pi n^{3/2} $
then Theorem~\ref{planar.curv}
implies that $w_2(\gamma) > n$.
\qed

\section{Width in higher dimensions and codimensions}  \label{general}

We now define width for a general submanifold.
We first consider a submanifold $\Gamma$ and
a collection of submanifolds transverse to $\Gamma$. 
We then consider a more general situation where $\Gamma$ is not
necessarily transverse to a family of submanifolds.

\noindent
{\bf Definition 1.} 
Let $N$ and $\Gamma$ be submanifolds of a manifold $M$ and let $G$ be a $k$-dimensional group of diffeomorphisms of $M$.   Let $X$ denote the set of
submanifolds $\{ g(N): g \in G  \} $ that transversely intersect $\Gamma$.
Define an equivalence relation on $X$  by setting $g_0(N) \sim g_1(N) $ if there is a path from  $g_0$ to $g_1$ in $G$ such that $g_t(N)  , 0 \le t \le 1$ is in $X$. So $g_0(N) \sim g_1(N) $ if  $g_0(N)$ is isotopic to $g_1(N)$ through submanifolds $g_t(N)$ transverse to $\Gamma$.
 
Equivalent submanifolds $g_0(N) $ and $ g_1(N) $  intersect $\Gamma$ in the
same number of points.  Pick one representative $g_i(N)$ for each equivalence class. The {\em width $w_G$} of  $\Gamma$ is defined to be the sum over equivalence classes of the number 
of intersections of $g_i(N)$ and $\Gamma$;
$$
w_G(\Gamma) = \sum_i  |\Gamma \cap g_i(N)|.
$$
 
 We say that $\Gamma$ is   {\em $w_G$-minimizing} if $\Gamma$ minimizes $w_G$ within
 its isotopy class.
 
 \noindent
{\bf Example 1.}  Let $M, N, G$ be $\RR^3$, the $xy$-plane and $\RR^1$ respectively. For $g \in \RR^1$ let $g ( (x,y,0)) = (x,y,g)$.  Then $G$ defines a foliation of $\RR^3$ by horizontal planes.  
A curve $\Gamma$ is $w_G$-minimizing if it is in thin position as in \cite{Gabai}.
 
 \noindent
 {\bf Example 2.} Let $M, N, G$ be ${\mathbb S}^3 \times {\mathbb S}^1$, ${\mathbb S}^3 \times \{1\}$,  and ${\mathbb S}^1$,  respectively, and let $ {\mathbb S}^1$ act on $M$ by rotation.  A curve $\gamma$ in
 $M$  is in thin position relative to $w_G$ if it either lies in a leaf neighborhood with a single maximum and minimum and has width two
 or is always transverse to the leaves. In the latter case there is a single equivalence class and the homotopy class of $\gamma$ determines its width, equal to  $\omega ( \gamma)$, where $\omega$ generates  $H^1(M;Z)$.
 
It is not necessary to restrict our attention to transverse intersections.  We can obtain still more general measures of intersection complexity by using the following definition.

\noindent
{\bf Definition 2.} 
Let $N$ and $\Gamma$ be submanifolds of a manifold $M$ and let $G$ be a $k$-dimensional group of diffeomorphisms acting on $M$.   Let $X$ denote the set of
submanifolds $\{ g(N) : g \in G  \} $  and
define an equivalence relation on $X$  by setting $g_0(N) \sim g_1(N) $ if there is a path from  $g_0$ to $g_1$ in $G$ such that $g_t(N) \cap  \Gamma$ is  isotopic to
$g_0(N) \cap  \Gamma, ~ 0 \le t \le 1$. So $g_0(N) \sim g_1(N) $ if  $g_0(N)$ is isotopic to $g_1(N)$ through submanifolds whose intersection with $\Gamma$ is preserved up to isotopy.  (Transversal intersection is no longer required.) The width of a submanifold $\Gamma$ is
the intersection number of $\Gamma$ and $\{ g(N):g \in G \}$ counted with appropriate multiplicity.  
 
 \noindent
{\bf Example 3.}   Let $M, N, G$ be $\RR^3$, the $z$-axis and $\RR^2$ respectively. For $g \in \RR^2$ let $g ((0,0,z)) = (g,z)$.  Then $G$ defines a foliation of $\RR^3$ by vertical lines.  Let $\gamma$
be a smooth curve. 

For a generic $\gamma$, the equivalence classes of $X_G$ consist of 
\begin{itemize}
\item Sets of lines that miss $\gamma$, one for each complementary component of the
projection $\pi(\gamma)$ to the $xy$-plane.
\item One component for each edge of   $\pi(\gamma)$.
\item One component for each vertex of   $\pi(\gamma)$.
\end{itemize}
Define the width of a line $g \cdot N$ in $X_G$ to be $n(n-1)$ if $g \cdot N \cap \gamma \ = n$.
The width here is chosen so that perturbing $\gamma$ to a generic projection which
has only double point singularities preserves the width. Only crossings contribute to the
width, which is equal to the crossing number of the projection of $\gamma$  to the
$xy$-plane.  The 2-width of a knot is equal to the
knot's crossing number.

 \noindent
{\bf Example 4.}   Let $M, N, G$ be $\RR^3$, the $xy$-plane and $\RR^1$ respectively
as in Example 1. But we now allow non-transverse intersections and take the width of 
$\{ g(N): g \in G  \} $ to be   the intersection number of $\Gamma$ and $\{ g(N): g \in G  \} $ counted without multiplicity.  
 
For a generic $\gamma$, the equivalence classes of $X_G$ consist of 
\begin{itemize}
\item Two components represented by planes that miss $\gamma$, one above and one below $\gamma$.
\item One component for each equivalence class intersecting $\gamma$ transversely.
\item One component for each plane intersecting $\gamma$ non-transversely at one or more points.
\end{itemize}

The width minimizing representative of an $n$-bridge knot is a curve $\gamma$ in bridge position, with all its maximums at one level
and all its minimum at a second level. A width-minimizing curve has $w_G = 4n$
where $n$ is the bridge number.

\begin{flushright} 

Joel Hass\\ Department of Mathematics\\ University
of California\\ Davis,
CA 95616\\ e-mail: hass@math.ucdavis.edu\\
\end{flushright}

\begin{flushright} 
J. Hyam Rubinstein\\ Department of Mathematics and Statistics\\
University of Melbourne\\ Parkville, 3010
Australia\\ e-mail: rubin@ms.unimelb.edu.au\\
\end{flushright}

\begin{flushright} 
Abigail Thompson\\ Department of Mathematics\\
University of California\\ Davis,
CA 95616\\ e-mail: thompson@math.ucdavis.edu\\
\end{flushright}

\end{document}